\documentclass[11pt,letterpaper,reqno]{amsart}

\usepackage{amsmath}
\usepackage{amsfonts}
\usepackage{amssymb}
\usepackage{amsthm}
\usepackage{amscd}

\usepackage[colorlinks=true]{hyperref}

\usepackage{latexsym}
\usepackage{mathrsfs}
\usepackage{psfrag}
\usepackage[dvips]{epsfig}
\usepackage{epsfig}
\usepackage[all]{xy}         

\usepackage[margin=3.0cm]{geometry}

\theoremstyle{plain}                              
\newtheorem{thm}{Theorem}[section]

\newtheorem{defn}[thm]{Definition}
\newtheorem{lem}[thm]{Lemma}
\newtheorem{cor}[thm]{Corollary}

\newtheorem{question}{Question}

\newtheorem*{oseledets}{Oseledets Multiplicative Ergodic Theorem}
\newtheorem*{thmA}{Theorem A}
\newtheorem*{thmB}{Theorem B}
\newtheorem*{thmC}{Theorem C}

\theoremstyle{definition}                         
\newtheorem*{example}{Example}
\newtheorem*{remark}{Remark} 

\theoremstyle{remark}                             

\numberwithin{equation}{section}

\newcommand{\R}{\mathbb{R}}                     

\newcommand{\ms}[1]{\mathscr{#1}}               

\newcommand{\veps}{\varepsilon}

\newcommand{\uvar}{\norm{u}_\infty-\chi}

\providecommand{\norm}[1]{\left\lVert#1\right\rVert}       

\providecommand{\abs}[1]{\left\lvert#1\right\rvert}        

\begin{document}

\title[Oseledets regularity functions]{Oseledets regularity functions
  for Anosov flows}


\author[S. N. Simi\'c]{Slobodan N. Simi\'c}

\address{Department of Mathematics, San Jos\'e State University, San
  Jos\'e, CA 95192-0103}



\email{simic@math.sjsu.edu}


\date{\today}
\dedicatory{}



\begin{abstract}
  Oseledets regularity functions quantify the deviation of the growth
  associated with a dynamical system along its Lyapunov bundles from
  the corresponding uniform exponential growth. Precise degree of
  regularity of these functions is unknown. We show that for every
  invariant Lyapunov bundle of a volume preserving Anosov flow on a
  closed smooth Riemannian manifold, the corresponding Oseledets
  regularity functions are in $L^p(m)$, for some $p > 0$, where $m$ is
  the probability measure defined by the volume form. We prove an
  analogous result for essentially bounded cocycles over volume
  preserving Anosov flows.

\end{abstract}

\maketitle





\section{Introduction}
\label{sec:intro}

This paper is concerned with the so called Oseledets regularity
functions, which naturally arise from the Oseledets Multiplicative
Ergodic Theorem and the Birkhoff Ergodic Theorem. We restrict our
attention to flows; a similar analysis could be done for
diffeomorphisms.

For a smooth flow $\Phi = \{ f_t \}$ on a closed Riemannian manifold
$M$ and a nonzero vector $v \in T_x M$, the \textsf{Lyapunov exponent}
of $v$ is defined by
\begin{displaymath}
  \chi(v) = \lim_{\abs{t} \to \infty} \frac{1}{t} \log \norm{T_x f_t(v)},
\end{displaymath}
if the limit exists. Vectors $v$ with the same Lyapunov exponent
$\chi$ (plus the zero vector) form a linear subspace $E^\chi(x)$ of
$T_x M$, called the Lyapunov space of $\chi$. By construction, these
spaces form an invariant bundle in the sense that $T_x f_t(E^\chi(x))
= E^\chi(f_t x)$, for all $t \in \R$.

The fundamental properties of Lyapunov spaces are described by the
following seminal result (originally stated for diffeomorphisms).

\begin{oseledets}[\cite{osel+68,ruelle+79,barreira+pesin+01}]
  Let $\Phi = \{ f_t \}$ be a $C^1$ flow on a closed Riemannian
  manifold $M$. There exists a $\Phi$-invariant set $\ms{R} \subset M$
  of full measure with respect to any $\Phi$-invariant Borel
  probability measure $\mu$, such that for every $x \in \ms{R}$ there
  exists a splitting (called the Oseledets splitting)
\begin{displaymath}
  T_x M = \bigoplus_{i=1}^{\ell(x)} E_i(x),
\end{displaymath}
and numbers $\chi_1(x) < \cdots < \chi_{\ell(x)}(x)$ with the
following properties:

\begin{enumerate}

\item[(a)] The bundles $E_i$ are $\Phi$-invariant,
  \begin{displaymath}
    T_x f_t(E_i(x)) = E_i(f_tx),
  \end{displaymath}
  and depend Borel measurably on $x$.

\item[(b)] For all $v \in E_i(x) \setminus \{0\}$, 
  \begin{displaymath}
    \lim_{\abs{t} \to \infty} \frac{1}{t} \log \norm{T_x f_t (v)} = \chi_i(x).
  \end{displaymath}
  The convergence is uniform on the unit sphere in $E_i(x)$.

\item[(c)] For for any $I, J \subset \{ 1, \ldots, \ell(x) \}$ with $I
  \cap J = \emptyset$, the angle function is tempered, i.e.,
  \begin{displaymath}
    \lim_{\abs{t} \to \infty} \frac{1}{t} \log 
    \sphericalangle(T_x f_t(E_I(x)),T_x f_t(E_J(x))) = 0,
  \end{displaymath}
  where $E_I = \bigoplus_{i \in I} E_i$.

\item[(d)] For every $x \in \ms{R}$, 
  \begin{displaymath}
    \lim_{\abs{t} \to \infty} \frac{1}{t} \log \det T_x f_t = \sum_{i=1}^{\ell(x)} 
      \chi_i(x) \dim E_i(x).
  \end{displaymath}
\item[(e)] If $\Phi$ is ergodic with respect to $\mu$, then the
  functions $\ell$ and $\chi_i$ are $\mu$-almost everywhere constant.

\end{enumerate}
\end{oseledets}

Points $x \in \ms{R}$ are called \textsf{regular}. 

Assume $\Phi$ is ergodic with respect to some measure $\mu$, fix $i
\in \{1, \ldots, \ell \}$, and set $\chi = \chi_i$ and $E =
E_i$. Denote the restriction of $Tf_t$ to $E$ by $T^E f_t$. Since
\begin{displaymath}
  \chi = \lim_{t \to \infty} \frac{1}{t} \log \norm{T_x^E f_t},
\end{displaymath}
for each $x \in \ms{R}$, it follows that for every $\veps > 0$,
\begin{displaymath}
  \lim_{t \to \infty} \frac{\norm{T_x^E f_t}}{e^{(\chi+\veps)t}} = 0.
\end{displaymath}
Therefore, there exists a constant $C > 0$, depending on $x$ and
$\veps$, such that $\norm{T_x^E f_t} \leq C e^{(\chi+\veps)t}$, for
all $t \geq 0$. It is natural to consider the best such $C$ as a
function of $x$ and $\veps$:

\begin{defn}
  For a fixed Lyapunov bundle $E$ of $\Phi$ and every $\veps > 0$, the
  $(E,\veps)$-\textsf{Oseledets regularity function} $R_\veps : \ms{R}
  \to \R$ is defined by
  \begin{displaymath}
    R_\veps(x) = \sup_{t \geq 0} 
            \frac{\norm{T_x^E f_t}}{e^{(\chi + \veps)t}}.
  \end{displaymath}
\end{defn}

The family $\{R_\veps\}_{\veps > 0}$ is the main focus of this
paper. It is not hard to see that each $R_\veps$ is Borel measurable
(see \cite{barreira+pesin+01} for the case of diffeomorphisms) and
that $R_\veps \geq 1$. What more can be said about the $R_\veps$? For
instance:

\begin{question}
  Does $R_\veps$ lie in some $L^p$-space? What is the best value of
  $p$?
\end{question}

A related question can be posed for cocycles. Recall that a map
$\Delta : M \times \R \to \R$ is called a (multiplicative real-valued)
cocycle over a flow $\{ f_t \}$ if
\begin{displaymath}
  \Delta(x,s+t) = \Delta(x,s) \Delta(f_s x, t),
\end{displaymath}
for all $x \in M$ and $s, t \in \R$. If for every $x \in M$ the map $t
\mapsto \Delta(x,t)$ is absolutely continuous, then
\begin{equation}     \label{eq:Delta}
  \Delta(x,t) = \exp \left\{ \int_0^t u(f_s x) \: ds \right\},
\end{equation}
where $u(x) = \dot{\Delta}(x,0) = \frac{d}{dt} \left|_{t=0} \log
  \Delta(x,t) \right.$. When $u$ is essentially bounded with respect
to some measure, we will call such a cocycle \textsf{essentially bounded}.

Assume $\mu$ is an invariant Borel probability measure, $u \in
L^\infty(\mu)$, and set $\chi = \int_M u \: d\mu$. If $\mu$ is
ergodic, then by the Birkhoff Ergodic Theorem,
\begin{displaymath}
  \lim_{t \to \infty} \frac{1}{t} \log \Delta(x,t) = \chi,
\end{displaymath}
for $\mu$-a.e. $x$. Denote the set of Birkhoff regular points (at
which the above limit exists and equals $\chi$) by $\ms{R}$ as
well. Then for every $x \in \ms{R}$ and $\veps > 0$, $\Delta(x,t)/
\exp\{ (\chi + \veps)t \} \to 0$, as $t \to \infty$, so there exists a
constant $B > 0$ (depending on $x$ and $\veps$) such that $\Delta(x,t)
\leq B \exp\{ (\chi + \veps)t \}$, for all $t \geq 0$. It makes sense
to consider the best such $B$ as a function of $x$ and $\veps$:

\begin{defn}
  For each $\veps > 0$, the $(u,\veps)$-\textsf{regularity function}
  $D_\veps^u : M \to \R$ is defined by
  \begin{equation}     \label{eq:D}
    D_\veps^u(x) = \sup_{t \geq 0} \frac{\Delta(x,t)}{e^{(\chi+\veps)t}}.
  \end{equation}
\end{defn}

When $u$ is clear from the context, we will write just $D_\veps$. It
is clear that each $D_\veps$ is Borel measurable and $D_\veps \geq
1$. It is also easy to see that if $\veps \geq \uvar$, then $D_\veps =
1$, $\mu$-a.e. We are therefore interested only in the values of
$\veps$ less than $\uvar$.

What more can be said about the $D_\veps$? For instance:

\begin{question}
  Does $D_\veps$ belong to some $L^p$-space? What is the best value of
  $p$?
\end{question}

A word of caution is in place here. Even for ``best'' (interesting)
dynamical systems, namely globally uniformly hyperbolic ones, neither
the Oseledets nor the Birkhoff theorem guarantee any particularly good
properties of the set $\ms{R}$ of regular points and as a consequence
of the regularity functions. Although $\ms{R}$ has full measure with
respect to any invariant probability measure, its complement is not
only non-empty, but can be topologically very large. This follows from
the work of Barriera and Schmeling~\cite{barr+schmel+00} who showed
that for an Anosov diffeomorphism of the 2-torus, the complement of
the set of regular points can have the full Hausdorff dimension (i.e.,
two). In the continuous time case, using multifractal analysis,
Barreira and Saussol~\cite{barreria+saussol+00} showed that for
hyperbolic flows the set of non-regular points is similarly
topologically large, namely, dense and of full Hausdorff dimension,
for a generic function $u$. Similar results were obtained by Pesin and
Sadovskaya~\cite{pesin+sad+01}.

We will soon see that Questions 1 and 2 are closely related, at least
in the case of Anosov flows (see \S~\ref{sec:anosov-flows}), to which
we now restrict ourselves. Namely, given a volume preserving Anosov
flow $\Phi$ and a Lyapunov bundle $E$ of $\Phi$, it turns out that
each regularity function $R_\veps$ of $E$ can be related to a
regularity function $D_\eta^u$, for some $\eta > 0$ and some
essentially bounded function $u$ dependent only on $E$. See Theorem B.

We now state our main results. Throughout, $m$ will denote the Borel
probability measure defined by the Riemannian volume on $M$.

\begin{thmA}
  Let $\Phi = \{ f_t \}$ be a $C^2$ volume preserving Anosov flow on a
  closed Riemannian manifold $M$ and let $\Delta : M \times \R \to \R$
  be a multiplicative cocycle over $\Phi$, as in \eqref{eq:Delta}. If
  $u \in L^\infty(m)$, then for every $\veps > 0$, the corresponding
  $(u,\veps)$-regularity function $D_\veps$ belongs to $L^p(m)$, for
  some $p > 0$.

  If $u$ is H\"older continuous, let $H$ be the entropy function of
  $u$ (as defined in \S~\ref{sec:anosov-flows}). Then $D_\veps \in
  L^p(m)$, provided that
  \begin{equation}    \label{eq:H-int}
    p \leq \left\{ \int_\veps^{\uvar} \frac{ds}{H(\chi+s)} \right\}^{-1}.
  \end{equation}
\end{thmA}

Here is a sketch of the proof of Theorem A. If $u$ is H\"older, then
for all $x \in \ms{R}$, $t \mapsto \Delta(x,t)$ is continuous, so we
define $T_\veps : \ms{R} \to \R$ ($0 < \veps < \uvar$) to be the
smallest $t \geq 0 $ at which the supremum in \eqref{eq:D} is
attained. Then $T_\veps$ is Borel measurable and $D_\veps \leq \exp \{
(\uvar - \veps) T_\veps\}$, so we study the question of integrability
of $\exp(T_\veps)$. Using a large deviations result of
Waddington~\cite{wadd+96} (see \S~\ref{sec:anosov-flows} for details),
we show that if $p < H(\chi + \veps)$, then $\exp(T_\veps) \in L^p(m)$,
where $H$ is the entropy function of $u$. Next, we show that if $\eta
< \veps$, then $D_\eta \leq D_\veps \exp \{ (\veps - \eta) T_\eta \}$
$m$-a.e., which for any natural number $N$ by induction extends to
\begin{displaymath}
  D_\veps \leq \prod_{i=0}^{N-1} \exp(\delta T_{\veps_i}),
\end{displaymath}
where $\veps = \veps_0 < \veps_1 < \cdots < \veps_N = \uvar$ is a
partition of $[\veps, \uvar]$ with $\delta = \veps_{i+1} - \veps_i =
(\uvar-\veps)/N$, for all $i$. Using the generalized H\"older
inequality and the fact that $\exp(\delta T_{\veps_i}) \in
L^{p_i}(m)$, where $p_i < H(\chi + \veps_i)/\delta$, we obtain
$D_\veps \in L^p(m)$, where $p^{-1} = \sum_i p_i^{-1} > \sum_i
\delta/H(\chi+\veps_i)$. Passing to the limit as $N \to \infty$ in the
last sum, we obtain \eqref{eq:H-int}.

If $u$ is only essentially bounded, we show that it is possible to
suitably approximate $u$ in the $L^1$-sense by a larger smooth
function (Lemma~\ref{lem:smooth}). Namely, for every $\delta > 0$
there exists a $C^\infty$ function $\tilde{u} : M \to \R$ such that $u
\leq \tilde{u}$ and $\int_M (\tilde{u} - u) \: dm < \delta$. It then
easily follows that for any $0 < \delta < \veps < \uvar$, $D_\veps^u
\leq D_{\veps-\delta}^{\tilde{u}}$, $m$-a.e., which implies that
$D_\veps^u$ lies in some $L^p$-space.

A bridge between the two different types of regularity functions is
given by the following result.

\begin{thmB}      
  Let $\Phi = \{ f_t \}$ be a $C^2$ volume preserving Anosov flow on a
  closed $C^\infty$ Riemannian manifold $M$ and let $E$ be a Lyapunov
  bundle for $\Phi$ associated with a Lyapunov exponent $\chi$. For
  every $\delta > 0$ there exists a constant $C_\delta > 0$ such that
  \begin{displaymath}
    \norm{T_x^E f_t} \leq C_\delta e^{\delta t} 
    \exp \left\{ \int_0^t u(f_s x) \: ds \right\},
  \end{displaymath}
  for every $x \in \ms{R}$ and $t \geq 0$, where $u \in L^\infty(m)$
  is independent of $\delta$ and 
  \begin{displaymath}
    \int_M u \: dm = \chi.
  \end{displaymath}
\end{thmB}

The proof of Theorem B goes as follows. First, we trivialize $E$ by
using a measurable orthonormal frame. This transforms the second
variational equation for the restriction of the flow to $E$ into a
family of non-autonomous differential equations $\dot{X} = A_x(t) X$
on $\R^k$ ($k = \dim E$), parametrized by $x \in \ms{R}$. Following
\cite{barreira+pesin+01}, we use a lemma of Perron to construct a
family $U_x(t)$ of orthogonal matrices such that if $v(t)$ is a
solution to $\dot{v} = A_x(t) v$, then $z(t) = U_x(t) v(t)$ is a
solution to $\dot{z} = B_x(t) z$, where $B_x(t)$ are upper triangular
matrices, whose non-diagonal entries are bounded in $x$ and $t$. We
show that for every $\delta > 0$ there exists a norm
$\norm{\cdot}_\delta$ on $\R^k$ such that $\norm{B_x(t)}_\delta <
r(B_x(t)) + \delta$, where $r$ denotes the spectral radius. Moreover,
$r(B_x(s+t)) = r(B_{f_s x}(t))$, for all $t$ and
a.e. $x$. Furthermore, if $X(t)$ is the unique solution to the matrix
differential equation $\dot{X} = A_x(t) X$ satisfying $X(0) = I$, then
\begin{displaymath}
  \norm{X(t)}_\delta \leq K_\delta e^{\delta t} \exp \left\{
    \int_0^t r(B_x(s)) \: ds \right\},
\end{displaymath}
for some constant $K_\delta > 0$. We therefore define $u : M \to \R$
by $u(x) = r(B_x(0))$. It is not hard to prove that $u$ is essentially
bounded. The desired inequality for $T_x^E f_t$ is now obtained by
pulling the norms $\norm{\cdot}_\delta$ back to $E$ and observing that
the each new Finsler structure is globally uniformly equivalent to the
original one.

The last result is a straightforward corollary of Theorem B.

\begin{thmC}
  Let $\Phi = \{ f_t \}$ be a $C^2$ volume preserving Anosov flow on a
  closed $C^\infty$ Riemannian manifold $M$. Let $E$ be a Lyapunov
  bundle in the Oseledets splitting for $\Phi$. Then for every $\veps
  > 0$, the corresponding $(E,\veps)$-regularity function $R_\veps$
  belongs to the space $L^p(m)$, for some $p > 0$.
\end{thmC}

To prove Theorem C, denote the Lyapunov exponent corresponding to $E$
by $\chi$ and let $\veps > 0$ and $0 < \delta < \veps$ be
arbitrary. Then by Theorem B,
\begin{displaymath}
  \frac{\norm{T_x^E f_t}}{e^{ (\chi + \veps) t }} 
    \leq C_\delta 
    \frac{ \exp \left\{ \int_0^t u(f_s x) \: ds \right\}}
      {e^{(\chi+\veps-\delta)t}} \\
       \leq C_\delta D_{\veps-\delta}(x),
\end{displaymath}
for $m$-a.e. $x \in \ms{R}$ and $t \geq 0$. This implies that $R_\veps
\leq C_\delta D_{\veps-\delta}$, which yields Theorem C.

\begin{remark}
  The question of the best $p = p(\veps)$ such that $D_\veps \in
  L^p(m)$ (and the analogous question for $R_\veps$) remains open. It
  is likely that the answer can be found by a more careful analysis of
  the set $\ms{L} = \{ (\veps, p) : D_\veps \in L^p(m) \}$, which
  possesses a number of interesting properties such as:

  \begin{itemize}
  \item[(a)] The set $\ms{L}' = \{ (\veps, p^{-1}) : (\veps, p) \in
    \ms{L} \}$ is convex. To see this, first observe that
    $D_{(1-\alpha) \veps_0 + \alpha \veps_1} \leq
    D_{\veps_0}^{1-\alpha} D_{\veps_1}^\alpha$, for all $\veps_0,
    \veps_1 > 0$ and $0 \leq \alpha \leq 1$ (the proof is
    straightforward). If $(\veps_i, p_i^{-1}) \in \ms{L}'$, $i = 1,
    2$, and $0 < \alpha < 1$, then $D_{\veps_i} \in L^{p_i}$ ($i =
    1,2$), so $D_{\veps_0}^{1-\alpha} \in L^{p_0/(1-\alpha)}$ and
    $D_{\veps_1} \in L^{p_1/\alpha}$. The above inequality and
    H\"older's inequality yield $D_\veps \in L^p$, where $\veps =
    (1-\alpha) \veps_0 + \alpha \veps_1$ and $p^{-1} = (1-\alpha)
    p_0^{-1} + \alpha p_1^{-1}$. Thus $\ms{L}'$ contains the line
    segment connecting $(\veps_0,p_0^{-1})$ and $(\veps_1,p_1^{-1})$.

  \item[(b)] If $\nu$ is a Borel probability measure on an interval $I
    \subset (0,\uvar)$ and $\phi : I \to \R$ a positive Borel function
    whose graph is contained in $\ms{L}$, then $(\veps,p) \in \ms{L}$,
    where $\veps = \int_I t \: d\nu(t)$ and $p = \left( \int_I
      d\nu/\phi \right)^{-1}$.
  \end{itemize}
\end{remark}

The following example shows that even for ``simple'' systems there is
a definite cut-off value for $p$ beyond which the regularity function
is not in $L^p$.

\begin{example}
  Let $M = T^2$ be the 2-torus and $f : T^2 \to T^2$ an area
  preserving Anosov diffeomorphism. Denote its a.e. Lyapunov exponents
  by $\chi^- < 0 < \chi^+$. It is not hard to construct $f$ so that it
  possesses two periodic points $x, y$, whose corresponding Lyapunov
  exponents are different, i.e., $\chi_x^- \neq \chi_y^-$ and
  $\chi_x^+ \neq \chi_y^+$ (with obvious notation). Clearly, the
  Lyapunov exponents at $x$ or $y$ (or both) differ from the
  a.e. Lyapunov exponents $\chi^-, \chi^+$. Assume, for instance, that
  $\chi_x^+ > \chi^+$ and denote by $\lambda^+$ the unstable cocycle
  of $f$, that is, the determinant of the derivative of $f$ restricted
  to the unstable bundle of $f$. Then for $0 < \veps < \chi^+_p -
  \chi^+$, the regularity function $D_\veps$ of $u = \log \lambda^+$
  is infinite at $x$. (Using the fact that the homoclinic points of
  $x$ are dense in $T^2$ and all have the same Lyapunov exponents as
  $x$, it is not hard to see that $D_\veps$ is in fact infinite on a
  dense subset of $T^2$.)

  We claim that there exists $p_0 > 0$ such that $D_\veps \not\in
  L^p(m)$, for all $p \geq p_0$. The main idea for showing this is the
  following. Since $\lambda^+$ is continuous, each $D_\veps$ is
  lower-semicontinuous, so sets $E_\alpha = \{ D_\veps > \alpha \}$
  are open, for all $\alpha$. Since $D_\veps(x) = \infty$ (where $x$
  is the periodic point as above), $x \in E_\alpha$, for all $\alpha$,
  so there exists a ball $B(x,r)$, for some $r = r(\alpha)$, contained
  in $E_\alpha$. H\"older continuity of $u$ allows us to control the
  size of $r$ and show that $\int_0^\infty \alpha^{p-1} m(E_\alpha) \:
  d\alpha$ diverges for large enough $p$. The details follow.

  Denote by $C > 0$ and $0 < \theta < 1$ the H\"older constant and
  exponent of $u$ so that for all $x_1, x_2 \in T^2$,
  \begin{displaymath}
    \abs{u(x_1) - u(x_2)} \leq C d(x_1,x_2)^\theta.
  \end{displaymath}
  Fix $0 < \veps < \chi_x^+ - \chi^+$ and define $\sigma = \chi_x^+ -
  \chi^+ - \veps$ and 
  \begin{displaymath}
    S_\veps(z,N) = \sum_{i=0}^{N-1} u(f^iz) - (\chi^+ + \veps)N,
  \end{displaymath}
  so that $D_\veps(z) = \sup_{N \geq 1} \exp S_\veps(z,N)$. Consider
  the periodic point $x$ as above, at which $D_\veps(x) =
  \infty$. Suppose its prime period is $\ell$. It is not hard to
  verify that $S_\veps(x,n\ell) = \sigma n \ell$, for all $n \geq 1$,
  that is, $S_\veps(x,N)$ grows linearly along the subsequence $N = n
  \ell$. Moreover, for all $z \in T^2$, we have
  \begin{equation}      \label{eq:S_eps}
    \abs{S_\veps(x,N) - S_\veps(z,N)}  \leq \sum_{i=0}^{N-1}
    \abs{u(f^ix) - u(f^iz)} 
    \leq C \frac{\lambda^{\theta N} - 1}{\lambda^\theta - 1} d(x,z)^\theta,
  \end{equation}
  where $\lambda > 1$ is the Lipschitz constant of $f$. For $\alpha >
  0$, set $E_\alpha = \{ z \in T^2 : D_\veps(z) > \alpha \}$, as
  above. Clearly, $x \in E_\alpha$ for all $\alpha > 0$. Since
  $S_\veps(x,n \ell) = \sigma n \ell$, it follows that
  $S_\veps(x,n\ell) > \log \alpha$, where we take
  \begin{displaymath}
    n = \left[ \frac{\log \alpha}{\sigma \ell} \right] + 1.
  \end{displaymath}
  Fix $N = n \ell$ and observe that
  \begin{equation}   \label{eq:N}
    \frac{\log \alpha}{\sigma} + \ell \leq N < \frac{\log \alpha}{\sigma} 
        + 2\ell.
  \end{equation}
  It easily follows from \eqref{eq:S_eps} that if the right-hand side
  in \eqref{eq:S_eps} is $< S_\veps(x,N) - \log \alpha$, then
  $S_\veps(z,N) > \log \alpha$, hence $z \in E_\alpha$. Thus the ball
  $B(x,r(\alpha))$ in $T^2$ of radius
  \begin{displaymath}
    r(\alpha) = \left\{ \frac{1}{C} 
      \frac{\lambda^\theta - 1}{\lambda^{N \theta} - 1} 
      [S_\veps(x,N) - \log \alpha] \right\}^{1/\theta}
  \end{displaymath}
  is contained in $E_\alpha$. It follows from \eqref{eq:N} that
  \begin{displaymath}
    \lambda^{N \theta} < \lambda^{\theta \left( \frac{\log \alpha}{\sigma}
        + 2 \ell \right)}  
     = \lambda^{2 \ell \theta} \alpha^{\frac{\theta \log \lambda}{\sigma}} 
     =: C_1 \alpha^{\rho},
  \end{displaymath}
  where $\rho = \theta \log \lambda/\sigma$. Similarly, $S_\veps(x,N)
  - \log \alpha = \sigma N - \log \alpha \geq \sigma \ell$. Hence
  \begin{displaymath}
    r(\alpha) \geq \left\{ \frac{1}{C} 
      \frac{\lambda^\theta - 1}{C_1 \alpha^{\rho} - 1} 
      \sigma \ell \right\}^{1/\theta} 
    \geq \left\{ \frac{1}{C} 
      \frac{\lambda^\theta - 1}{C_1 \alpha^{\rho}} 
      \sigma \ell \right\}^{1/\theta}
    = C_2 \alpha^{- \log \lambda/\sigma},
  \end{displaymath}
  and thus $m(E_\alpha) \geq m(B(x,r(\alpha)) \geq \pi C_2^2 \alpha^{-
    2\log \lambda/\sigma}$. Since
  \begin{displaymath}
    \int_{T^2} D_\veps^p \: dm = \int_0^\infty \alpha^{p-1} m(E_\alpha) \: d\alpha,
  \end{displaymath}
  we conclude that $D_\veps \not\in L^p(m)$, if $p \geq \frac{\log
    \lambda}{\sigma} = (\log \lambda)/(\chi_x^+ - \chi^+ - \veps)$. \qed
\end{example}

The paper is organized as follows. In Section~\ref{sec:prelim} we
recall some basics facts about Anosov flows, present the large
deviation result of Waddington~\cite{wadd+96}, and review some
Pesin-Lyapunov theory used later in the paper. The proofs of Theorems
A and B are given in Sections~\ref{sec:proofA} and \ref{sec:proofB}.

\subsection*{Submultiplicative cocycles}
\label{sec:sub}

A similar argument can be extended to all (essentially bounded)
\emph{submultiplicative} cocycles over Anosov flows, that is, maps $A
: M \times \R \to \R$ such that
\begin{displaymath}
  A(x,s+t) \leq A(x,s) A(f_s x, t),
\end{displaymath}
for all $x \in M$ and $s, t \in \R$. By Kingman's subadditive ergodic
theorem \cite{king+68} applied to $a = \log A$, it follows that for
a.e. $x$,
\begin{displaymath}
  \lim_{t \to \infty} \frac{1}{t} a(x,t) = \chi = \inf_{t > 0}
  \frac{1}{t} \int_M a(x,t) \: dm(x).
\end{displaymath}
If $\chi$ is finite, we can define regularity functions of $A$ as
above by
\begin{displaymath}
  R_\veps(x) = \sup_{t > 0} \frac{A(x,t)}{e^{(\chi + \veps) t}}.
\end{displaymath}
The goal is to show that for every $\veps > 0$, $R_\veps \in L^p(m)$,
for some $p > 0$. We briefly outline how this could be done. 

The key is to obtain the asymptotics of $m \{ R_\veps > e^\alpha \}$
with respect to $\alpha$. Choose $T_0 > 0$ large enough so that
\begin{displaymath}
  \frac{1}{T_0} \int_M a(x,T_0) \: dm(x) < \chi + \frac{\veps}{2}.
\end{displaymath}
If $R_\veps(x) > e^\alpha$, then $a(x,T) - (\chi + \veps) T > \alpha$,
for some $T > 0$. Write $T = kT_0 + \tau$, for some positive integer
$k$ and $0 \leq \tau < T_0$ and observe that $a(x,T)$ is bounded above
by the sum of $a(f_{T_0}^k x, \tau)$ (which is bounded a.e.) and the
$k^{th}$ Birkhoff sum of the function $g(x) = a(x,T_0)$. Thus the set
$\{ R_\veps > e^\alpha \}$ is contained in the set of the form $\{
\sum_{i=0}^{k-1} g \circ f_{T_0}^i > c \}$, for some $c$ depending on
$\alpha$, so one can apply to $f_{T_0}$ the large deviations result
for time-$t$ maps of Anosov flows proved by Dolgopyat in
\cite{dolgopyat+04}, yielding the desired asymptotics.

This approach (which we will not pursue here) could be used to prove
both Theorems A and C without the Pesin-Lyapunov theory, although it
does not provide the more precise bound on $p$ as a function of
$\veps$ given in Theorem A in the H\"older case. We thank the referee
for pointing out the possibility of this extension.

\section{Preliminaries}
\label{sec:prelim}

\subsection{Large deviations for  Anosov flows}
\label{sec:anosov-flows}

A non-singular $C^1$ flow $\Phi = \{ f_t \}$ on a closed Riemannian
manifold $M$ is called an \textsf{Anosov flow} if there exists a
$Tf_t$-invariant continuous splitting of the tangent bundle,
\begin{equation*}
  TM = E^{uu} \oplus E^c \oplus E^{ss},
\end{equation*}
and constants $C, \lambda > 0$ such that for all $t \geq 0$,
\begin{equation*}
\norm{Tf_t \! \restriction_{E^{ss}}} \leq C e^{-\lambda t} 
\qquad \qquad \text{and} \qquad
\qquad \norm{Tf_{-t} \! \restriction_{E^{uu}}} \leq C e^{-\lambda t},
\end{equation*}
where the center bundle $E^c$ is one dimensional and generated by the
infinitesimal generator of the flow. The bundles $E^{uu}, E^{ss}$ are
called the \textsf{strong unstable} and \textsf{strong stable bundle}
of the flow. If the flow is of class $C^2$, $E^{ss}, E^{uu}$ are known
to be H\"older continuous (cf., \cite{hassel+94, hassel+97,
  hps77}). If an Anosov flow preserves a volume form, it is
automatically ergodic with respect to the Lebesgue measure defined by
the volume (see \cite{anosov+67}). Recall that a flow is called (topologically)
transitive if it has a dense orbit. 

An \textsf{equilibrium state} of a function $\varphi : M \to \R$ is an
invariant Borel probability measure $\mu$ at which the quantity
\begin{displaymath}
  h(\mu) + \int_M \varphi \: d\mu
\end{displaymath}
attains its supremum, where $h(\mu)$ denotes the measure-theoretic
entropy of $\Phi$ with respect to $\mu$. This supremum $P(\varphi)$
is called the \textsf{pressure} of $\varphi$. If $\varphi$ is H\"older
continuous, there exists a \emph{unique} equilibrium state of
$\varphi$, denoted by $\mu_\varphi$.

Given a transitive Anosov flow $\Phi = \{ f_t \}$, one defines a
function $\varphi^u : M \to \R$ by
\begin{displaymath}
  \varphi^u(x) = \left. \frac{d}{dt} \right|_0 \log \det T_x f_t \! \! 
     \restriction_{E^{uu}}.
\end{displaymath}
If $\Phi$ is $C^2$, $\varphi^u$ is known to be H\"older
continuous. The unique equilibrium state of $-\varphi^u$ is called the
\textsf{Sinai-Ruelle-Bowen (SRB)} measure $\mu_\text{\footnotesize
  SRB}$ of the flow. By the Bowen-Ruelle
theorem~\cite{bowen-ruelle+75}, for every continuous $\varphi : M \to
\R$,
\begin{displaymath}
  \lim_{T \to \infty} \frac{1}{T} \int_0^T \varphi(f_t x) \: dt 
     = \int_M \varphi \: d\mu_\text{\rm SRB},
\end{displaymath}
for $m$-a.e. $x \in M$, where $m$ is the Lebesgue measure defined by
the volume form and $\Phi$ is $C^2$. Thus $\mu_\text{\footnotesize
  SRB}$ is an ergodic measure for $\Phi$. If the flow admits a
smooth invariant Borel probability measure $\mu$ (i.e., a measure that
is absolutely continuous with respect to the volume measure $m$), then
by the Birkhoff ergodic theorem, $\mu = \mu_\text{\footnotesize
  SRB}$. In particular, if $\Phi$ is volume preserving, then
$\mu_\text{\footnotesize SRB} = m$.

For an arbitrary flow $f_t : M \to M$ and function $\psi : M \to \R$,
we can define a skew product flow
\begin{displaymath}
  S_t^\psi : S^1 \times M \to S^1 \times M
\end{displaymath}
by
  \begin{displaymath}
    S_t^\psi(\exp (2\pi i \theta), x)
    = \left( \exp\{ 2\pi i (\theta + \psi^t(x)) \}, f_t(x) \right),
  \end{displaymath}
  where
  \begin{displaymath}
    \psi^t(x) = \int_0^t \psi(f_s x) \: ds.
  \end{displaymath}

\begin{defn}[\cite{wadd+96}]
  A H\"older continuous function $\varphi : M \to \R$ and a flow
  $\Phi = \{ f_t \}$ on $M$ are called \textsf{flow independent} if
  they satisfy the following property: for every two numbers $a, b \in
  \R$, if the skew product flow $S_t^{a+b\varphi}$ is \emph{not}
  topologically transitive\footnote{Waddington uses the term
    topologically ergodic, which has the same meaning, see
    \cite{peter+89}, Proposition 2.4.}, then $a = 0 = b$.
\end{defn}

Large deviation asymptotics for transitive Anosov flows were
established by Waddington in \cite{wadd+96}. In particular:

\begin{thm}[Corollary 2, \cite{wadd+96}]   \label{thm:wadd}
  Let $\Phi = \{ f_t \}$ be a transitive $C^2$ Anosov flow on $M$ and
  let $\varphi : M \to \R$ be a H\"older continuous function such that
  $\varphi$ and $\Phi$ are flow independent. There exist analytic
  real-valued functions $\beta, \gamma, \rho$ defined on an interval
  in $\R$, such that if $\rho(a) > 0$, then
  \begin{displaymath}
    \mu_\text{\rm SRB} 
    \left\{ x : \int_0^T \varphi(f_t x) \: dt 
      \geq Ta \right\}
      \sim \frac{C(a)}{\rho(a)} \frac{1}{2\pi \beta''(\rho(a))} 
      \frac{e^{\gamma(a) T}}{\sqrt{T}},
  \end{displaymath}
  as $T \to \infty$, where $C(a)$ is a constant depending on $a$.
\end{thm}

Here, $a(t) \sim b(t)$ as $t \to \infty$, means $a(t)/b(t) \to 1$. The
function $\beta : \R \to \R$ is defined by $\beta(t) = P(\psi + t
\varphi) - P(\psi)$, for a H\"older continuous $\psi$. For our
purposes, we will take $\psi = 0$.

Some properties of $\beta$ (see \cite{wadd+96} for details), with
$\psi = 0$, are:
\begin{equation}    \label{eq:beta}
  \beta'(t) = \int_M \varphi \: d\mu_{t\varphi} \quad
  \text{and} \quad \beta''(t) = \sigma^2_{\mu_{t\varphi}}(\varphi),
\end{equation}
where for a measure $\mu$ with $\int_M \varphi \: d\mu = \chi$, the
variance of $\varphi$ is defined by
\begin{displaymath}
  \sigma^2_\mu(\varphi) = \lim_{T \to \infty} \frac{1}{T} \left(
    \int_0^T \varphi \circ f_t \: dt - \chi T \right)^2.
\end{displaymath}
Furthermore, $\sigma^2_\mu(\varphi) = 0$ if and only if $\varphi$ is
cohomologous to a constant. If $\varphi$ is not cohomologous to a
constant, the map $t \mapsto \beta'(t)$ is strictly increasing. Denote
its range by $\Gamma_\varphi$; it follows from \eqref{eq:beta} that
$\Gamma_\varphi \subset (\min \varphi, \max \varphi)$. On
$\Gamma_\varphi$, set $\rho = (\beta')^{-1}$. Then $\rho :
\Gamma_\varphi \to \R$ is strictly increasing, surjective, and real
analytic, with $\rho(\chi) = 0$. Finally, $\gamma : \Gamma_\varphi \to
\R$ is defined as minus one times the Legendre transform of $\beta$,
i.e.,
\begin{displaymath}
  \gamma(s) = - \sup_{t \in \R} \{st - \beta(t) \}.
\end{displaymath}
It can be shown that $\gamma$ is a strictly concave, non-positive
function with a unique maximum at $\chi = \int_M \varphi \: dm$
(where we still take $\psi = 0$). Furthermore, $\gamma''(s) =
-1/\beta''(\rho(s))$, so in particular, $\gamma''(\chi) = -
1/\sigma^2_m(\varphi)$ (cf., \cite{wadd+96}).

In the large deviations literature the function $H = -\gamma$ is
called the \textsf{entropy function} of $\varphi$. It is easily seen
that $H$ has the following properties (see \cite{wadd+96}): it is
strictly convex on $\Gamma_\varphi$,
\begin{displaymath}
  H(\chi) = H'(\chi) = 0, \quad H''(\chi) = \frac{1}{\sigma_m^2(\varphi)},
  \quad \text{and} \quad 
  H(a) = \infty \ \text{for} \ a \not\in \Gamma_\varphi,
\end{displaymath}
where $\chi = \int_M \varphi \: dm$.

The following lemma will be needed later in the paper.

\begin{lem}      \label{lem:flow_ind}
  Let $\Phi$ be a volume preserving Anosov flow and $\varphi : M \to \R$ a
  H\"older continuous function. If $\varphi$ and $\Phi$ are not flow
  independent, then $\varphi$ is cohomologous to a constant.
\end{lem}

\begin{proof}
  Suppose $\varphi$ and $\Phi$ are not flow independent. Then there
  exist numbers $a, b$, not both zero, such that the skew product
  $S_t^{a + b \varphi}$ is not topologically transitive, hence not
  ergodic with respect to the measure $m_1 \times m$, where $m_1$ is
  the Haar-Lebesgue measure on $S^1$. Since the volume measure is an
  equilibrium state of $\Phi$, Proposition 4.2 in \cite{walkden+99}
  implies the existence of a nonzero integer $\ell$ and a H\"older
  function $w : M \to \R$ such that
  \begin{displaymath}
    \ell \int_0^t (a + b \varphi)(f_s x) = w(f_t x) - w(x),
  \end{displaymath}
  for all $x \in M$. If $b = 0$, then $a \neq 0$ and $w(f_t x) - w(x)
  = a \ell t$ everywhere, which is impossible. Therefore, $b \neq
  0$. Differentiating the above identity, we obtain
  \begin{displaymath}
    \varphi + \frac{a}{b} = \frac{1}{\ell} Xw,
  \end{displaymath}
  which means that $\varphi$ is cohomologous to $-a/b$. 
\end{proof}


\subsection{Pesin-Lyapunov theory}
\label{sec:lyapunov-theory}

In this section we follow Barreira-Pesin~\cite{barreira+pesin+01} and
briefly review some elements of Pesin-Lyapunov theory for linear
differential equations
\begin{equation}          \label{eq:ode}
  \dot{v} = A(t) v,
\end{equation}
where $A(t)$ is a $k \times k$ bounded matrix function, i.e.,
\begin{displaymath}
  \sup_{t \in \R} \norm{A(t)} < \infty.
\end{displaymath}
We concentrate on \emph{real} matrices $A(t)$
(\cite{barreira+pesin+01} deals with complex matrices). The Lyapunov
exponent of $v \in \R^k$ is the number
\begin{displaymath}
  \chi(v) = \limsup_{t \to \infty} \frac{1}{t} \log \norm{v(t)},
\end{displaymath}
where $v(t)$ is the unique solution to \eqref{eq:ode} satisfying the
initial condition $v(0) = v$. The function $\chi : \R^k \to \R \cup \{
- \infty \}$ attains only finitely many values $\chi_1 < \ldots <
\chi_\ell$, where $\ell \leq k$. For each $1 \leq i \leq \ell$, define
\begin{displaymath}
  V_i = \{ v \in \R^k : \chi(v) \leq \chi_i \}.
\end{displaymath}
This defines a linear filtration of $\R^k$:
\begin{displaymath}
  \{ \mathbf{0} \} = V_0 \subsetneqq V_1 \subsetneqq \cdots 
      \subsetneqq V_\ell = \R^k.
\end{displaymath}
An ordered basis $\mathbf{v} = (v_1, \ldots, v_k)$ of $\R^k$ is called
\textsf{normal} with respect to the filtration $\ms{V} = \{ V_i \}$ if
for every $1 \leq i \leq \ell$, the vectors $v_1, \ldots, v_{k_i}$
form a basis for $V_i$, where $k_i = \dim V_i$. In particular, if
$\chi$ is a constant function, every basis of $\R^k$ is normal.

Given a basis $\mathbf{v} = (v_1, \ldots, v_k)$ of $\R^k$ and $1 \leq
m \leq k$, denote by $\Gamma_m^\mathbf{v}(t)$ the volume of the
parallelepiped defined by $v_1(t), \ldots, v_m(t)$, where $v_i(t)$ is
the unique solution to \eqref{eq:ode} satisfying $v_i(0) =
v_i$. Recall that the Lyapunov exponent $\chi$ is regular (together
with the Lyapunov exponent $\tilde{\chi}$ associated with the dual
equation $\dot{w} = - A(t)^\ast w$) if and only if (see
\cite{barreira+pesin+01}, Theorem 1.3.1) for any normal ordered basis
$\mathbf{v} = (v_1, \ldots, v_k)$ of $\R^k$ and $1 \leq m \leq k$, we
have
\begin{displaymath}
  \lim_{t \to \infty} \frac{1}{t} \log \Gamma_m^\mathbf{v}(t) 
     = \sum_{i=1}^m \chi(v_i).
\end{displaymath}
In particular, if $\chi$ is constant, then for any basis $\mathbf{v} =
(v_1, \ldots, v_k)$ and $1 \leq m \leq k$,
\begin{equation}    \label{eq:m-chi}
  \lim_{t \to \infty} \frac{1}{t} \log \Gamma_m^\mathbf{v}(t) = m \chi.
\end{equation}
We now recall how one converts (as in \cite{barreira+pesin+01}) by a
linear change of coordinates the equation \eqref{eq:ode} into $\dot{z}
= B(t) z$, where the matrix $B(t)$ is \emph{upper triangular}. We seek
a differentiable family of orthogonal matrices $U(t)$ for the job. Set
$z(t) = U(t)^{-1} v(t)$, where $v(t)$ is a solution to \eqref{eq:ode};
then
\begin{displaymath}
  \dot{v}(t) = \dot{U}(t) z(t) + U(t) \dot{z}(t) = A(t) U(t) z(t),
\end{displaymath}
which implies $\dot{z}(t) = B(t) z(t)$, where
\begin{equation}      \label{eq:B(t)}
  B(t) = U(t)^{-1} A(t) U(t) - U(t)^{-1} \dot{U}(t).
\end{equation}
The following lemma of Perron guarantees the existence of $U(t)$ so
that $B(t)$ is upper triangular, for all $t$.

\begin{lem}[Lemma 1.3.3, \cite{barreira+pesin+01}]  \label{lem:perron}
  There exists a differentiable matrix function $t \mapsto U(t)$ such
  that for each $t \geq 0$:
  \begin{enumerate}

  \item[(a)] $U(t)$ is orthogonal.

  \item[(b)] $B(t) = [b_{ij}(t)]$ is upper triangular.

  \item[(c)] For all $1 \leq i < j \leq k$, 
    \begin{displaymath}
      \sup_{t \geq 0} \abs{b_{ij}(t)} \leq 2 \sup_{t \geq 0} \norm{A(t)} < \infty.
    \end{displaymath}
  \item[(d)] For any basis $\mathbf{v} = (v_1, \ldots, v_k)$ of $\R^k$
    and all $1 \leq i \leq k$,
    \begin{displaymath}
      b_{ii}(t) = \frac{d}{dt} \log \frac{\Gamma_i^\mathbf{v}(t)}
          {\Gamma_{i-1}^\mathbf{v}(t)}.
    \end{displaymath}
  \end{enumerate}
\end{lem}

Here is how families $U(t)$ and $B(t)$ are constructed. Denote by 
\begin{displaymath}
  \ms{G} : Gl(k,\R) \to O(k)
\end{displaymath}
the Gram-Schmidt orthogonalization operator that sends a basis
$\mathbf{v} = (v_1,\ldots,v_k)$ of $\R^k$ to an orthonormal basis
$\mathbf{u} = (u_1,\ldots,u_k)$. We can think of $\mathbf{v}$ and
$\mathbf{u}$ as matrices with columns $v_1,\ldots,v_k$ and
$u_1,\ldots, u_k$, respectively. Then $\mathbf{v} \in Gl(k,\R)$ and
$\mathbf{u} \in O(k)$. Observe that 
\begin{displaymath}
  \ms{G} = \ms{N} \circ \ms{L},
\end{displaymath}
where $\ms{L}[v_1,\ldots,v_k] = [w_1,\ldots,w_k]$ is a linear operator
defined by
\begin{displaymath}
  w_{i+1} = v_{i+1} - \text{proj}_{W_i} v_{i+1}, \qquad 
     W_i = \text{span}\{ w_1, \ldots, w_i\},
\end{displaymath}
and $\ms{N}[w_1,\ldots,w_k] = [u_1,\ldots, u_k]$ is the normalization
operator
\begin{displaymath}
  u_i = \frac{w_i}{\norm{w_i}}.
\end{displaymath}
Since $\ms{L}$ is linear, differentiating $\ms{G}$ at $v$ yields
\begin{displaymath}
  T_{\mathbf{v}} \ms{G} = T_{\ms{L}\mathbf{v}} \ms{N} \circ \ms{L}.
\end{displaymath}
In the proof of Perron's Lemma 1.3.1 in \cite{barreira+pesin+01}, for
an arbitrary but fixed basis $\mathbf{v} = (v_1,\ldots,v_k)$ of
$\R^k$, one defines $U(t)$ by
\begin{displaymath}
  U(t) = \ms{G}[v_1(t),\ldots, v_k(t)],
\end{displaymath}
where $v_i(t)$ is the unique solution to the equation $\dot{v} = A(t)
v$ satisfying the initial condition $v(0) = v_i$.

The family $B(t)$ is defined as in \eqref{eq:B(t)}. Thus both $t
\mapsto U(t)$ and $t \mapsto B(t)$ depend on the choice of a basis
$\mathbf{v} = (v_1,\ldots,v_k)$ of $\R^k$. When it is important to
emphasize this, we will write $U^\mathbf{v}(t)$ and $B^\mathbf{v}(t)$.

It is clear that the eigenvalues of $B(t)$ are its diagonal entries
$b_{ii}(t)$. Denote the spectral radius of a matrix $M$ by $r(M)$.

\begin{cor}        \label{cor:r}
  If $\chi$ is constant, then
  \begin{displaymath}
    \lim_{t \to \infty} \frac{1}{t} \int_0^t r(B(s)) \: ds = \chi.
  \end{displaymath}
\end{cor}

\begin{proof}
  Follows directly from \eqref{eq:m-chi} and part (d) of
  Lemma~\ref{lem:perron}.
\end{proof}

\begin{lem}      \label{lem:independence}
  The spectral radius $r(B(t))$ of the matrix $B(t) =
  B^{\mathbf{v}}(t)$ is independent of the choice of the basis
  $\mathbf{v} = (v_1,\ldots, v_k)$ of $\R^k$.
\end{lem}

\begin{proof}
  Let $\mathbf{v} = (v_1,\ldots, v_k)$ and $\mathbf{w} = (w_1,\ldots,
  w_k)$ be two bases of $\R^k$ and let $B^{\mathbf{v}}(t)$ and
  $B^{\mathbf{w}}(t)$ be the corresponding matrices constructed as
  above. Denote the solutions to \eqref{eq:ode} with initial values
  $(v_1,\ldots, v_k)$, $(w_1,\ldots, w_k)$ by $(v_1(t),\ldots,v_k(t))$
  and $(w_1(t), \ldots,w_k(t))$, respectively. Both $k$-tuples are
  bases of $\R^k$. Therefore, there exists a family of invertible
  matrices $P(t)$ such that $P(t) v_i(t) = w_i(t)$, for all $1 \leq i
  \leq k$. It follows that
  \begin{displaymath}
    \Gamma_i^{\mathbf{w}}(t) = \det P(t) \: \Gamma_i^{\mathbf{v}}(t),
  \end{displaymath}
  for all $1 \leq i \leq k$, and thus
  \begin{displaymath}
     \frac{\Gamma_i^\mathbf{w}(t)}
          {\Gamma_{i-1}^\mathbf{w}(t)} = 
           \frac{\Gamma_i^\mathbf{v}(t)}
          {\Gamma_{i-1}^\mathbf{v}(t)},
  \end{displaymath}
  for all $t \geq 0$. Lemma~\ref{lem:perron}(d) implies that the
  corresponding diagonal entries of $B^{\mathbf{v}}(t)$ and
  $B^{\mathbf{w}}(t)$ are the same, which yields the conclusion of the
  lemma.
\end{proof}

Define a function $\rho_B : \R \to \R$ by
\begin{displaymath}
  \rho_B(t) = r(B(t)).
\end{displaymath}

\begin{lem}       \label{lem:rho}
  There exists a universal constant $K > 0$, depending only on $k$,
  such that
  \begin{displaymath}
    \abs{\rho_B(0)} \leq K \norm{A(0)}.
  \end{displaymath}
\end{lem}

\begin{proof}
  Let $\mathbf{v} = \mathbf{e}$ be the standard basis
  $(e_1,\ldots,e_k)$ of $\R^k$ and let $U(t) = U^{\mathbf{e}}(t)$ be
  the corresponding orthogonal matrix function defined as above. Then:
  \begin{align*}
    \abs{\rho_B(0)} & = \abs{r(B(0))} \\
      & = \abs{r(A(0) - U(0)^{-1} \dot{U}(0))} \\
      & \leq \norm{A(0) -  U(0)^{-1} \dot{U}(0)} \\
      & \leq \norm{A(0)} + \norm{U(0)^{-1} \dot{U}(0)} \\
      & = \norm{A(0)} + \norm{\dot{U}(0)}.
  \end{align*}
  Denote the solution to \eqref{eq:ode} with initial value $e_i$ by
  $e_i(t)$. Then:
  \begin{align*}
    \dot{U}(0) & = \left. \frac{d}{dt} \right|_0 
       \ms{G}[e_1(t),\ldots,e_k(t)] \\
       & = T_I \ms{G} [\dot{e}_1(0), \ldots, \dot{e}_k(0)] \\
       & = T_I \ms{G} [A(0) e_1, \ldots, A(0) e_k] \\
       & = T_I \ms{G} (A(0)),
  \end{align*}
  where $I$ is the $k \times k$ identity matrix. Let $K = 1 + \norm{T_I
    \ms{G}}$, where $T_I \ms{G}$ is regarded as a map between Lie
  algebras $\mathfrak{gl}_k$ and $\mathfrak{o}_k$. It follows that
  \begin{displaymath}
    \rho_B(0) \leq K \norm{A(0)},
  \end{displaymath}
  completing the proof of the lemma.
\end{proof}

\section{Proof of Theorem A}
\label{sec:proofA}

We split the proof of Theorem A into two cases: in the first case, we
deal with H\"older continuous functions $u$. The general case of
essentially bounded $u$ is reduced to the first case in a suitable
way. In either case, without loss of generality, we assume that $u$
\textit{is a positive function}. Otherwise, apply the analysis below
to the function $u + C$, for a sufficiently large positive constant
$C$. It is easy to see that the regularity functions of $u$ and $u+C$
are the same. \\

\noindent
\textbf{Case 1:} \textit{$u$ is H\"older continuous.} \\

First of all, we may assume that $u$ and $\Phi$ are flow
independent. Otherwise, by Lemma~\ref{lem:flow_ind}, $u$ is
cohomologous to a constant, which is necessarily equal to $\chi =
\int_M u \: dm$, that is, $u = Xw + \chi$, for some H\"older function
$w$. This implies that
\begin{displaymath}
  \exp \left\{ \int_0^t u(f_s x) \: ds \right\} 
  = e^{\chi t} e^{w(f_t x) - w(x)} \leq e^{2 \norm{w}_\infty} e^{\chi t},
\end{displaymath}
so the corresponding regularity functions $D_\veps$ are all constant
(in fact, $D_\veps = 1$ $m$-a.e., for all $\veps > 0$). 

Recall that we are only interested in the values $0 < \veps <
\norm{u}_\infty - \chi$, since $D_\veps = 1$ for $\veps \geq \uvar$.

Denote the set of Birkhoff regular points by $\ms{R}$. For each $x \in
\ms{R}$ and $0 < \veps < \uvar$, define $T_\veps(x)$ to be the
smallest $T \geq 0$ at which the supremum defining $D_\veps =
D_\veps^u$ in \eqref{eq:D} is attained. That is, 
\begin{displaymath}
  T_\veps(x) = \min \left\{ T \geq 0 : \int_0^T u(f_s x) \: ds - (\chi + \veps)T
      = \log D_\veps(x) \right\}.
\end{displaymath}
By the Birkhoff Ergodic Theorem, $T_\veps : \ms{R} \to [0, \infty)$ is
well-defined. It is clear that $T_\veps$ is Borel measurable.

As in \S~\ref{sec:anosov-flows}, let $H = -\gamma$ be the entropy
function of $u$.

\begin{lem}      \label{lem:T}
  If $p < H(\chi + \veps)$, then $e^{T_\veps} \in L^p(m)$.
\end{lem}
\begin{proof}
  Fix an $\veps \in (0, \norm{u}_\infty - \chi)$ and let $\zeta > 1$
  be arbitrary. Define
  \begin{displaymath}
    B_n = \{ x : \zeta^n < T_\veps(x) \leq \zeta^{n+1} \}.
  \end{displaymath}
  Suppose $x \in B_n$. Since $u$ is assumed to be positive, we have
  \begin{align*}
    \int_0^{\zeta^{n+1}} u(f_s x) \: ds & \geq
    \int_0^{T\veps(x)} u(f_s x) \: ds \\
    & = (\chi + \veps) T_\veps(x) + \log D_\veps(x) \\
    & \geq (\chi + \veps) T_\veps(x) \\
    & \geq (\chi + \veps) \zeta^n \\
    & = \frac{\chi + \veps}{\zeta} \zeta^{n+1}.
  \end{align*}
  By Theorem~\ref{thm:wadd} there exists a constant $L$ depending on
  $\veps$ and $\zeta$ such that
  \begin{displaymath}
    m(B_n) \leq L \exp 
    \left\{ -H\left(\frac{\chi+\veps}{\zeta} \right)
      \zeta^{n+1} \right\}.
  \end{displaymath}
    It follows that 
    \begin{align*}
      \int_M \exp (p T_\veps) \: dm & = \sum_{n=0}^\infty \int_{B_n}
      \exp (p T_\veps) \: dm \\
      & \leq \sum_n L \exp \left\{ p \zeta^{n+1} 
        -H\left(\frac{\chi+\veps}{\zeta} \right)
        \zeta^{n+1} \right\},
    \end{align*}
    which is finite for 
    \begin{displaymath}
      p < H\left(\frac{\chi+\veps}{\zeta} \right). 
    \end{displaymath}
    Since $\zeta > 1$ was arbitrary, letting $\zeta \to 1+$ yields the claim.
  \end{proof}

  \begin{lem}       \label{lem:down}     
    If $0 < \eta < \veps$ and $x \in \ms{R}$, then
    \begin{displaymath}
      D_\eta(x) \leq D_\veps(x) e^{(\veps - \eta) T_\eta(x)}.
    \end{displaymath}
  \end{lem}
  \begin{proof}
    Set $u_\veps = u - \chi -\veps$. Then for each $x \in \ms{R}$:
    \begin{align*}
      \log D_\veps(x) & = \max_{t \geq 0} \int_0^t u_\veps(f_s x) \: ds \\
          & \geq \int_0^{T_\eta(x)} u_\veps(f_s x) \: ds \\
          & = \int_0^{T_\eta(x)} u_\eta(f_s x) \: ds + (\eta - \veps) T_\eta(x) \\
          & = \log D_\eta(x) - (\veps - \eta) T_\eta(x),
    \end{align*}
    which proves the claim.
  \end{proof}

  Now let $0 < \veps < \uvar$ be arbitrary and fix a natural number $N
  \geq 1$. Let
\begin{displaymath}
  \veps = \veps_0 < \veps_1 < \cdots < \veps_N = \uvar
\end{displaymath}
be a partition of the interval $[\veps, \uvar]$ with $\delta =
\veps_{i+1} - \veps_i = (\uvar - \veps)/N$, for all $0 \leq i \leq
N-1$. Applying Lemma~\ref{lem:down} repeatedly and using $D_{\uvar} =
1$ a.e., we obtain
\begin{displaymath}
  D_\veps \leq \prod_{i=0}^{N-1} \exp (\delta T_{\veps_i}).
\end{displaymath}
Since $\exp(\delta T_{\veps_i}) \in L^{p_i}(m)$, for $p_i <
H(\chi+\veps_i)/\delta$ (Lemma~\ref{lem:T}), the generalized H\"older
inequality yields $D_\veps \in L^p(m)$, where
\begin{align*}
  \frac{1}{p} & = \sum_{i=0}^{N-1} \frac{1}{p_i} \\
     & > \sum_{i=0}^{N-1} \frac{\delta}{H(\chi + \veps_i)} \\
     & \to \int_\veps^{\uvar} \frac{ds}{H(\chi+s)},
\end{align*}
as $N \to \infty$. This proves the second conclusion of Theorem A.

\begin{remark}
  It is possible to show that, in fact,
  \begin{displaymath}
    D_\veps(x) = \exp \left\{ \int_\veps^{\uvar} T_\eta(x) \: d\eta \right\},
  \end{displaymath}
  for $m$-almost every $x \in M$.
\end{remark}

\noindent
\textbf{Case 2:} $u \in L^\infty(m)$. \\

We need to show that for every $\veps > 0$ there exists $p > 0$ such
that $D_\veps \in L^p(m)$. The following lemma asserts that we may as
well work with a smooth $u$.

    \begin{lem}      \label{lem:smooth}
      For every $\delta > 0$ there exists a $C^\infty$ function
      $\tilde{u} : M \to \R$ such that $u \leq \tilde{u}$ and
    \begin{displaymath}
      \int_M (\tilde{u} - u) \: dm < \delta.
    \end{displaymath}
    \end{lem}

    \begin{proof}
      We will first find a continuous function $w : M \to \R$ such
      that $u \leq w$ and $\int (w-u) < \delta$ and then regularize
      $w$.

      Let $\eta > 0$ be arbitrary. By Luzin's theorem, there exists a
      continuous function $g : M \to [0, \infty)$ such that
      $\norm{g}_\infty \leq \norm{u}_\infty$ and $m(A) < \eta$, where
      $A = \{ x \in M : u(x) \neq g(x) \}$. The set $A$ is Borel
      measurable, so there exists an open set $U$ such that $A \subset
      U$ and $m(U \setminus A) < \eta$. Let $V$ be an open set such
      that
      \begin{displaymath}
        A \subset V \subset \overline{V} \subset U.
      \end{displaymath}
      By Urysohn's lemma, there exists a continuous function $h : M
      \to [0,1]$ such that $h = 0$ on the complement of $U$ and $h =
      1$ on $\overline{V}$. Let $k \in (\norm{u}_\infty, 2
      \norm{u}_\infty)$ and define
      \begin{displaymath}
        w = g + k h.
      \end{displaymath}
      On $A$, $w \geq k h = k > u$. On the complement of $A$, $g = u$, so
      $w = u + k h \geq u$. Observe that $w = u$ on the complement of
      $U$ and that $g + k h \leq 3 \norm{u}_\infty$ on $U$. Therefore,
      \begin{align*}
        \int_M (w - u) \: dm & = \int_U (g + k h - u) \: dm  \\
          & \leq \int_U (g+kh) \: dm  \\
          & \leq  m(U) \cdot 3 \norm{u}_\infty  \\
          & \leq 6 \eta \norm{u}_\infty.
      \end{align*}
      Let $w_a = w + a$, where $a > 0$ is a small constant, so that
      $w_a - u \geq a$. Finally, let $\tilde{u}$ be a $C^\infty$
      regularization of $w_a$ with $\norm{\tilde{u} - w_a}_\infty$
      sufficiently small so that $\tilde{u} \geq u$. It is easy to see
      that the integrals of $\tilde{u}$ and $w_a$ are the same. Since
      \begin{displaymath}
        \int_M (\tilde{u} - u) \: dm = \int_M (w_a - u) \: dm 
        \leq 6 \eta \norm{u}_\infty + a,
      \end{displaymath}
      by choosing $\eta$ and $a$ sufficiently small, we obtain a
      desired function $\tilde{u}$.
    \end{proof}

    Now fix an $\veps \in (0, \uvar)$ and $0 < \delta < \veps$. Let
    $\tilde{u}$ be a $C^\infty$ function on $M$ supplied by
    Lemma~\ref{lem:smooth} such that $u \leq \tilde{u}$ and
    $\tilde{\chi} - \chi < \delta$, where $\tilde{\chi} = \int
    \tilde{u} \: dm$. Denote by $\tilde{D}_\eta$ the
    $(\tilde{u},\eta)$-regularity function. Then:
    \begin{align*}
      D_\veps(x) & \leq \sup_{t \geq 0} 
        \frac{\exp \int_0^t \tilde{u}(f_s x) \: ds}{e^{(\chi + \veps) t}} \\
        & \leq \sup_{t \geq 0} 
           \frac{\exp \int_0^t \tilde{u}(f_s x) \: ds}
           {e^{(\tilde{\chi} + \veps - \delta) t}} \\
        & = \tilde{D}_{\veps - \delta}(x),
    \end{align*}
    which lies in $L^p(m)$, for some $p > 0$, by Case 1. Therefore $D_\veps
    \in L^p(m)$, completing the proof of Theorem A.

\section{Proof of Theorem B}
\label{sec:proofB}

Let $\Phi = \{ f_t \}$ be a $C^2$ volume preserving Anosov flow. Fix a
Lyapunov bundle $E$ corresponding to a Lyapunov exponent $\chi$ and
denote the set of Lyapunov regular points by $\ms{R}$.

Let $x \in \ms{R}$ and consider the Second Variational Equation for
the flow $\Phi$ on $E$:
  \begin{equation}   \label{eq:var}
    \frac{d}{dt} T_x^E f_t = (T_{f_t x}^E X) T_x^E f_t.
  \end{equation}
  where $X$ is the Anosov vector field. Choose a measurable
  orthonormal frame $\mathbf{F} = \{ F_1, \ldots, F_k \}$ for $E$ and
  define a vector bundle map
  \begin{displaymath}
    \ms{T} : E \to \ms{R} \times \R^k
  \end{displaymath}
  by
  \begin{displaymath}
    \ms{T}(F_i(x)) = (x, e_i),
  \end{displaymath}
  where $e_i$ is the $i^{th}$ element of the standard basis of $\R^k$;
  extend $\ms{T}$ linearly over each fiber. Then $\ms{T}$ trivializes
  $E$, transforming \eqref{eq:var} into a family of differential
  equations parametrized by $x \in \ms{R}$:
  \begin{displaymath}
    \dot{X} = A_x(t) X,
  \end{displaymath}
  where $A_x(t)$ is the matrix of $T_{f_t x}^E X$ relative to the
  frame $\mathbf{F}$. As in \S~\ref{sec:lyapunov-theory}, for
  each $x \in \ms{R}$ we obtain an orthogonal matrix $U_x(t)$ and an
  upper triangular matrix $B_x(t)$ whose properties are described by
  Lemma~\ref{lem:perron}. Observe that since
  \begin{displaymath}
    \sup_{x \in \ms{R}}
    \norm{T_x^E X} \leq \sup_{x \in M} \norm{T_x X} < \infty,
  \end{displaymath}
  it follows that
  \begin{equation}      \label{eq:A}
    \alpha = \sup \{ \norm{A_x(t)} : x \in \ms{R}, \: t \in \R \} < \infty.
  \end{equation}
  Furthermore, by Corollary~\ref{cor:r},
  \begin{displaymath}
     \lim_{t \to \infty} \frac{1}{t} \int_0^t r(B_x(s)) \: ds = \chi.
  \end{displaymath}
  It is well-known that for every matrix $M$ and $\delta > 0$ there
  exists a norm such that $\norm{M} < r(M) + \delta$. The following
  lemma is a slight generalization of this result.
  
  \begin{lem}    \label{lem:norm}
    Let $\beta > 0$ be fixed and denote by $\ms{B}$ the set of all
    upper triangular $k \times k$ matrices such that for all $B =
    [b_{ij}] \in \ms{B}$,
    \begin{displaymath}
      \max_{i < j} \abs{b_{ij}} \leq \beta.
    \end{displaymath}
    Then for every $\delta > 0$ there exists a norm $\norm{ \cdot
    }_\delta$ on $\R^k$ such that for all $B \in \ms{B}$, the induced
    operator norm of $B$ satisfies
    \begin{displaymath}
      \norm{B}_\delta < r(B) + \delta.
    \end{displaymath}
  \end{lem}

  \begin{proof}
    The proof is an adaptation of one of the standard proofs (see,
    e.g., \cite{kress+98}, Theorem 3.32). Define
    \begin{displaymath}
      \veps = \min \left\{ 1, \frac{\delta}{(k-1) \beta} \right\}
    \end{displaymath}
    and
    \begin{displaymath}
      D = \text{diag}(1,\veps, \veps^2, \ldots, \veps^{k-1}).
    \end{displaymath}
    Then for any $B = [b_{ij}] \in \ms{B}$,
    \begin{displaymath}
      C = D^{-1} B D =
      \begin{bmatrix}
        b_{11} & \veps b_{12} & \veps^2 b_{13} & \cdots 
        & \veps^{k-1} b_{1k} \\
        0 & b_{22} & \veps b_{23} & \cdots & \veps^{k-2} b_{2k} \\
        0 & 0 & b_{33} & \cdots & \veps^{k-3} b_{3k} \\
        \cdot & \cdot & \cdot & \cdots & \cdot \\
        0 & 0 & 0 & \cdots & b_{kk}
      \end{bmatrix}.
    \end{displaymath}
    For a $k \times k$ matrix $A =[a_{ij}]$, write
    \begin{displaymath}
      \norm{A}_\infty = \max_{1 \leq i \leq k} \sum_{j=1}^k \abs{a_{ij}}.
    \end{displaymath}
    Then, for all $B = [b_{ij}] \in \ms{B}$:
    \begin{displaymath}
      \norm{C}_\infty \leq \max_{1 \leq i \leq k} \abs{b_{ii}} 
        + (k-1) \veps \beta \leq r(B) + \delta.
    \end{displaymath}
    We define a norm on $\R^k$ by
    \begin{displaymath}
      \norm{v}_\delta = \norm{D^{-1} v}_\infty,
    \end{displaymath}
    where $\norm{(w_1,\ldots, w_k)}_\infty = \max \abs{w_i}$. It
    follows that for all $B \in \ms{B}$,
    \begin{align*}
      \norm{B v}_\delta & = \norm{D^{-1} B v}_\infty \\
        & = \norm{C D^{-1} v}_\infty \\
        & \leq \norm{C}_\infty \norm{D^{-1} v}_\infty \\
        & = \norm{C}_\infty \norm{v}_\delta \\
        & \leq (r(B) + \delta) \norm{v}_\delta. \qedhere
    \end{align*}
  \end{proof}

  By Lemma~\ref{lem:perron} and \eqref{eq:A}, 
  \begin{displaymath}
    \norm{B_x(t)} \leq 2 \norm{A_x(t)} \leq 2 \alpha,
  \end{displaymath}
  for all $x \in \ms{R}$ and $t \in \R$. Thus we can apply
  Lemma~\ref{lem:norm} to the family of matrices $\ms{B} = \{B_x(t) :
  x \in \ms{R}, t \in \R \}$. For each $\delta > 0$, we obtain a norm
  $\norm{ \cdot }_\delta$ on $\R^k$, which induces an operator matrix
  norm we also denote by $\norm{ \cdot }_\delta$. This yields
  \begin{displaymath}
    \norm{B_x(t)}_\delta \leq r(B_x(t)) + \delta,
  \end{displaymath}
  for all $t \in \R$ and $x \in \ms{R}$.

  Now consider the unique solution $X(t)$ to $\dot{X} = A_x(t) X$
  satisfying the initial condition $X(0) = I$ and the corresponding
  solution $Z(t) = U_x(t)^{-1} X(t)$ to $\dot{Z} = B_x(t) Z$. Since
  $\mathbf{F}$ is orthonormal, $U_x(0) = I$, so $Z(0) = I$. Thus
  \begin{displaymath}
    Z(t) = I + \int_0^t B_x(s) Z(s) \: ds.
  \end{displaymath}
  It follows that
  \begin{displaymath}
    \norm{Z(t)}_\delta \leq 1 + \int_0^t 
    \norm{B_x(s)}_\delta \norm{Z(s)}_\delta \: ds,
  \end{displaymath}
  so by Gr\"onwall's inequality and Lemma~\ref{lem:norm},
  \begin{equation}    \label{eq:Z}
    \norm{Z(t)}_\delta \leq \exp \left\{ \int_0^t \norm{B_x(s)}_\delta \: ds
        \right\} \leq e^{\delta t} \left\{ \int_0^t r(B_x(s)) \: ds
        \right\}.
  \end{equation}
  Since $U_x(t)$ is orthogonal, its operator norm with respect to the
  original norm on $\R^k$ equals one. The old norm and the new norm on
  $\R^k$ are uniformly equivalent, so there exists a uniform constant
  $K_\delta > 0$ such that $\norm{U_x(t)}_\delta \leq K_\delta
  \norm{U_x(t)} = K_\delta$. Therefore,
  \begin{equation}    \label{eq:X}
    \norm{X(t)}_\delta = \norm{U_x(t) Z(t)}_\delta 
         \leq K_\delta \norm{Z(t)}_\delta.
  \end{equation}
  Now let $\varrho(x,t) = r(B_x(t))$, for any choice of the matrices
  $B_x(t)$ as above. Note that by Lemma~\ref{lem:independence}
  $\varrho$ is well-defined.

\begin{lem}   \label{lem:flow_u}
  For all $x \in \ms{R}$ and $s, t \in \R$, we have
  \begin{displaymath}
    \varrho(x,s+t) = \varrho(f_s x,t).
  \end{displaymath}
\end{lem}

\begin{proof}
  Fix $x \in \ms{R}$. Recall how the matrices $B_x(t)$ are constructed
  (cf., \ref{eq:B(t)}): first, we choose a basis $\mathbf{v} =
  (v_1,\ldots,v_k)$ of $\{ x \} \times \R^k$, and apply the
  Gram-Schmidt procedure to the matrix $[v_1(t), \ldots,v_k(t)]$,
  where $\dot{v}_i(t) = A_x(t) v_i(t)$ and $v_i(0) = v_i$, which
  yields a family of orthogonal matrices $U_x^{\mathbf{v}}(t)$. Then
  we define $B_x^{\mathbf{v}}(t) = \{ U_x^{\mathbf{v}}(t) \}^{-1}
  A_x(t) U_x^{\mathbf{v}}(t) - \{ U_x^{\mathbf{v}}(t) \}^{-1}
  \dot{U}_x^{\mathbf{v}}(t)$.

  Now fix $s \in \R$. We define suitable families of matrices
  $B_x^{\mathbf{v}}(t)$ and $B_{f_s x}^{\mathbf{w}}(t)$ by
  appropriately choosing bases $\mathbf{v} = (v_1,\ldots,v_k)$ of
  $\ms{T}(E(x)) = \{x \} \times \R^k$ and $\mathbf{w} = (w_1, \ldots,
  w_k)$ of $\ms{T}(E(f_s x)) = \{ f_s x\} \times \R^k$,
  respectively. This can be done as follows.

  Define $\mathbf{v}$ by $(x,v_i) = \ms{T}(F_i(x))$ ($1 \leq i \leq
  k$). This gives rise to a family of orthogonal matrices
  $U_x^{\mathbf{v}}(t)$ and the corresponding family
  $B_x^{\mathbf{v}}(t)$.

  Define $\mathbf{w}$ by $(f_s x,w_i) = \ms{T}(T_x f_s (F_i(x)))$ ($1
  \leq i \leq k$). This gives rise to the matrices $U_{f_s
    x}^{\mathbf{w}}(t)$ and the corresponding family $B_{f_s
    x}^{\mathbf{w}}(t)$.

  Let $v_i(t)$ and $w_i(t)$ be the solutions of the differential
  equations $\dot{v} = A_x(t) v$ and $\dot{w} = A_{f_s x}(t) w$ with
  initial conditions $v_i$ and $w_i$, respectively. Then:
    \begin{align*}
      w_i(t) & = \ms{T}(T_{f_s x} f_t (T_x f_s (F_i(x))))  \\
      & = \ms{T}(T_x f_{s+t} (F_i(x))) \\
      & = v_i(s+t).
    \end{align*}
    This implies that $U_x^{\mathbf{v}}(s+t) = U_{f_s
      x}^{\mathbf{w}}(t)$.  Furthermore, since $A_x(t)$ is the matrix
    of $T_{f_t x}^E X$ in the frame $\mathbf{F}$, $A_x(t) = [ T_{f_t
      x}^E X]_{\mathbf{F}}$, it follows that
    \begin{displaymath}
      A_x(s+t) = [ T_{f_{s+t}x}^E X]_\mathbf{F} =  [ T_{f_t(f_s x)}^E X]_\mathbf{F}
      = A_{f_s x}(t).
    \end{displaymath}
    Therefore,
    \begin{align*}
      B_{f_s x}^{\mathbf{w}}(t) & = U_{f_s x}^{\mathbf{w}}(t)^{-1}
      A_{f_s x}(t) U_{f_s x}^\mathbf{w}(t) - U_{f_s x}^{\mathbf{w}}(t)
      \dot{U}^{\mathbf{w}}_{f_s x}(t) \\
      & = \left\{ U_x^{\mathbf{v}}(s+t) \right\}^{-1} A_x(s+t)
      U_x^{\mathbf{v}}(s+t)
      - \left\{ U_x^{\mathbf{v}}(s+t) \right\}^{-1} \dot{U}_x^{\mathbf{v}}(s+t) \\
      & = B_x^{\mathbf{v}}(s+t).
    \end{align*}
    It follows that $\varrho(x,s+t) = r(B_x^{\mathbf{v}}(s+t)) =
    r(B_{f_s x}^{\mathbf{w}}(t)) = \varrho(f_s x, t)$, as claimed.
\end{proof}

  Define a function $u : \ms{R} \to \R$ by
  \begin{displaymath}
    u(x) = r(B_x(0)) = \varrho(x,0).
  \end{displaymath}
  By Lemma~\ref{lem:rho},
  \begin{displaymath}
    \abs{u(x)} \leq K \norm{A_x(0)} \leq K \alpha,
  \end{displaymath}
  for $m$-a.e. $x$, hence $u \in L^\infty(m)$. Furthermore, by
  Lemma~\ref{lem:flow_u},
  \begin{equation}      \label{eq:u}
    u(f_t x) = r(B_{f_t x}(0)) = r(B_x(t)).
  \end{equation}
  Combining \eqref{eq:Z}, \eqref{eq:X}, and \eqref{eq:u}, we obtain
  \begin{displaymath}
    \norm{X(t)}_\delta \leq K_\delta e^{\delta t} \exp \left\{
         \int_0^t u(f_s x) \: ds \right\}.
  \end{displaymath}
  For each $\delta > 0$, we abuse the notation and denote the pullback
  of the norms $\norm{ \cdot }_\delta$ to $E$ via $\ms{T}$ by the same
  symbol. That is, for each $v \in E(x)$ ($x \in \ms{R}$), we set
  \begin{displaymath}
    \norm{v}_{\delta} = \norm{\ms{T}(v)}_\delta.
  \end{displaymath}
  This defines a family of Finsler structures on $E$ with respect to which 
  \begin{displaymath}
    \norm{T_x^E f_t}_\delta \leq K_\delta e^{\delta t} 
    \exp \left\{ \int_0^t u(f_s x) \: ds \right\},
  \end{displaymath}
  for all $x \in \ms{R}$ and $t \geq 0$.

  Since any two norms on $\R^k$ are uniformly equivalent, for each
  $\delta > 0$ there exists a constant $A_\delta > 0$ such that
  \begin{displaymath}
    \norm{v} \leq A_\delta \norm{v}_\delta,
  \end{displaymath}
  for all $v \in E$, where $\norm{v}$ denotes the original norm of $v$
  defined by the Riemann structure on $M$. It follows that the norm of
  $T_x^E f_t$ with respect to the original norm on $E$ satisfies
  \begin{displaymath}
    \norm{T_x^E f_t} \leq A_\delta \norm{T_x^E f_t}_\delta 
    \leq A_\delta K_\delta e^{\delta t} 
    \exp \left\{ \int_0^t u(f_s x) \: ds \right\}.
  \end{displaymath}
  This completes the proof of Theorem B.

\bibliographystyle{amsalpha} 

\end{document}